\documentclass[12pt,a4paper]{article}
\usepackage[a4paper, mag=1000,left=2.5cm, right=2.5cm, top=2cm, bottom=2cm, headsep=0.7cm, footskip=1cm,twoside]{geometry}
\usepackage{amsmath,amsfonts,amssymb,amsthm}
\usepackage{mathrsfs}
\usepackage{indentfirst,fancyhdr,ifthen,textcomp,multirow,longtable}
\usepackage[colorlinks=true,citecolor=magenta,linktocpage=true,linkcolor=blue,final]{hyperref}

\newtheorem{Theorem}{Theorem}[section]

\newtheorem{Def}{Definition}

\newtheorem{Remark}[Theorem]{Remark}

\newcommand\emp\varnothing
\newcommand\eps\varepsilon
\newcommand\ov\widetilde

\DeclareMathOperator{\pdiv}{div_1}
\DeclareMathOperator{\qdiv}{div_2}

\DeclareMathOperator{\pcurl}{curl}

\def\^#1{^{\overline{#1}}}
\newcommand{\email}[1]{\href{mailto:#1}{#1}}

\fancypagestyle{plain}{
	\fancyhead{}
	
	}
\newcommand{\Addresses}{{
  \bigskip
  \footnotesize

  O.\,P. Khromova,\par\nopagebreak \textsc{Altai State University}\par\nopagebreak \textsc{Barnaul, Russia}\par\nopagebreak
  \textit{E-mail address:} \email{khromova.olesya@gmail.com}

  \medskip

  P.\,N. Klepikov,\par\nopagebreak \textsc{Altai State University}\par\nopagebreak \textsc{Barnaul, Russia}\par\nopagebreak
  \textit{E-mail address:} \email{klepikov.math@gmail.com}

  \medskip

  S.\,V. Klepikova,\par\nopagebreak \textsc{Altai State University}\par\nopagebreak \textsc{Barnaul, Russia}\par\nopagebreak
  \textit{E-mail address:} \email{klepikova.svetlana.math@gmail.com}

  \medskip

  E.\,D. Rodionov,\par\nopagebreak \textsc{Altai State University}\par\nopagebreak \textsc{Barnaul, Russia}\par\nopagebreak
  \textit{E-mail address:} \email{edr2002@mail.ru}

}}

\begin{document}

\title{About the Schouten\nobreakdash--Weyl tensor\\on 3-dimensional Lorenzian Lie groups%
\thanks{This work was supported by the Russian Foundation for Basic Research (project \textnumero~16--01--00336a.)}}

\author{O.\,P. Khromova, P.\,N. Klepikov,\\ S.\,V. Klepikova, E.\,D. Rodionov}

\date{}

\maketitle

\begin{abstract}
The main purpose of this paper is to investigate the Schouten\nobreakdash--Weyl tensor on the~three-dimensional Lie groups with left-invariant Lorenzian metrics. The~left\nobreakdash-in\-vari\-ant Lorentzian metrics on the three\nobreakdash-dimensional Lie groups with squared length zero Schouten\nobreakdash--Weyl tensor are studied. Moreover, the three-dimensional metric Lie groups with almost harmonic (i.e.~with~zero curl and divergence) Schouten\nobreakdash--Weyl tensor are investigated. In addition, the~question about the harmonicity of contraction of~the~Schouten\nobreakdash--Weyl tensor is considered.
\end{abstract}

\allowdisplaybreaks

\section{Introduction}

\sloppy

The fundamental monography~\cite{Besse} give us a survey about the Riemannian metrics with~harmonic Weyl tensor. In dimension three, the Weyl tensor is trivial, therefore we consider the Schouten\nobreakdash--Weyl tensor (the Cotton tensor), which plays a role of~the~Weyl tensor. The Schouten\nobreakdash--Weyl tensor was investigated in~\cite{R-S-Sh} for the case of~left\nobreakdash-in\-vari\-ant Lorentzian metrics on the three-dimensional Lie groups. The previous paper was a~continuation of J.~Milnor's paper~\cite{Milnor} about the left\nobreakdash-invariant Riemannian metrics on~three\nobreakdash-dimensional Lie~groups. A~classification of the three-dimensional metric Lie~algebras of the Lie groups with a left-invariant Riemannian metric and trivial divergence and curl of the Schouten\nobreakdash--Weyl tensor is given in~\cite{G-R-S}.



\section{Basic notation and facts}

Let $G$ be a Lie group, $ \left\{ \mathfrak{g} ,[\cdot.\cdot] \right\}$ is a corresponding Lie algebra. We denote a~left\nobreakdash-in\-vari\-ant Lorentzian metric on $G$ by $\langle\cdot,\cdot\rangle$, Levi-Civita connection~by $\nabla$, curvature tensor and Ricci tensor by $R$ and $r$, which defined by the following
\begin{equation*}
R\left(X,Y\right)Z=\left[\nabla_Y,\nabla_X\right]Z+\nabla_{\left[X,Y\right]}Z, \quad 
r\left(X,Y\right)=\mathrm{tr}\left(Z\mapsto R\left(X,Z\right)Y\right).
\end{equation*}

Further we fix a basis $\{ E_{1},E_{2},\ldots ,E_{n}\}$ of left-invariant vector fields in $\mathfrak{g}$ and set
\begin{equation}\label{st_1}
[E_{i},E_{j}]=c_{ij}^{k}E_{k}, \quad \langle E_{i},E_{j}\rangle = g_{ij},
\end{equation}
where $\{c_{ij}^{k}\}$ are a structure constants of the Lie algebra and $\{g_{ij}\}$ is a metric tensor. 


We consider the one-dimensional curvature tensor $A$ and Schouten\nobreakdash--Weyl tensor $SW$, which defined by the formulas
\begin{equation*}\label{odn-kr}
A_{ik}=\frac{1}{n-2}\left( r_{ik}-\frac{\rho g_{ik}}{2(n-1)} \right), \quad SW_{ijk}=A_{ij,k}-A_{ik,j},
\end{equation*}
where $\rho$ is a scalar curvature, $A_{ij,k}$ are a covariant derivatives of $A_{ij}$.


The squared length of the Schouten\nobreakdash--Weyl tensor is defined by 
$$\|SW\|^2=SW_{ijk}SW^{ijk}.$$

\begin{Def} \label{def:0}
A tensor $T$ is isotropic,~if $\|T\|^2=0$ and $T\ne0$.
\end{Def}

\begin{Def} \label{def:1}
A tensor $T_{i_{1}\ldots i_{p}}$ 
is harmonic~(cf.~\cite[p.43]{Yano}),~if:
\begin{enumerate}
\item $T_{i_{1} \cdots  i_{p}}$  is antisymmetric,
\item $\mathrm{curl}(T_{ i_{1}i_{2} \cdots i_{p}})=0$ or      $T_{i_{1}i_{2} \cdots i_{p},t}=T_{t i_{2}\cdots i_{p},i_{1}} + T_{ i_{1}t      \cdots i_{p},i_{2}}+\cdots+T_{i_{1} i_{2}\cdots     t,i_{p}}$,
\item $\mathrm{div}(T_{ i_{1}i_{2} \cdots i_{p}})=g^{i_{1}t} T_{ i_{1} \ldots i_{p} ,t }=0$.
\end{enumerate}
\end{Def}

We note that the antisymmetric part $SW_{[ijk]}$ and symmetric part $SW_{(ijk)}$ of the~Schouten\nobreakdash--Weyl tensor vanish. Thus, condition 1) of Definition~\ref{def:1} is not satisfied by the~Schouten\nobreakdash--Weyl tensor.

\begin{Def} \label{def:2}
A tensor $T_{i_{1}\ldots i_{p}}$ 
is said to be almost harmonic~if:
\begin{enumerate}
\item $\mathrm{curl}(T_{ i_{1}i_{2} \cdots i_{p}})=0$,
\item $\mathrm{div}(T_{ i_{1}i_{2} \cdots i_{p}})=0$.
\end{enumerate}
\end{Def}

For the Schouten\nobreakdash--Weyl tensor $SW_{ijk}$ we introduce the divergence of type I and II by
\begin{equation*}
\pdiv (SW)=g^{it}SW_{ijk,t},\quad \text{and}\quad
\qdiv (SW)=g^{jt}SW_{ijk,t}.
\end{equation*}

Let $V=V^{i}E_{i}$ be a left-invariant vector field, which is identified with the vector $\{V^{k}\}$. Let $w_{ij}$ be the contraction of the Schouten\nobreakdash--Weyl tensor $SW_{kij}$ with the vector field $\{V^{k}\}$, i.e.
\begin{equation}\label{ev0}
w_{ij}=V^{k}SW_{kij}.
\end{equation}
Since the Schouten\nobreakdash--Weyl tensor $SW_{kij}$ is antisymmetric with respect to $i$ and~$j$, the tensor $w_{ij}$ is antisymmetric.
The covariant derivatives of $w_{ij}$ have the form
\begin{equation*}
w_{ij,k}=w_{lj}\Gamma^{l}_{ki}+w_{il}\Gamma^{l}_{kj}.
\end{equation*}
The curl and divergence of the tensor $w_{ij}$ are expressed as
\begin{equation*}
\begin{gathered}
\mathrm{curl}(w)=w_{ij,t}-w_{tj,i}-w_{it,j},\\
\mathrm{div}(w)=g^{it}w_{ij,t}.
\end{gathered}
\end{equation*}
The  length of the vector field $\{V^{k}\}$ squared is expressed as
\begin{equation}\label{lev1}
\|V\|^{2}=g_{ij}V^{i}V^{j},
\end{equation}
where $g_{ij}$ is the metric tensor of Lorentzian signature.

Together with the arbitrary vector fields we consider the harmonic vector fields.

\begin{Def} \label{def:3}
A vector field $\{V^{i}\}$ is called {\it harmonic} if:
\begin{enumerate}
\item $\mathrm{curl}(V)=V^{i}_{\hspace{4pt},j}-V^{j}_{\hspace{5pt},i}=0$ ,
\item $\mathrm{div}(V)=V^{i}_{\hspace{4pt},i}=0$,
\end{enumerate}
where the covariant derivatives of $\{V^{i}\}$ are found by
\begin{equation*}
V^{i}_{\hspace{4pt},k}=-V^{l}\Gamma^{i}_{lk}.
\end{equation*}
\end{Def}

Further classification results for three-dimensional Lorentzian Lie groups was obtained~in~\cite{R-S-Sh,PKh-2015}. 

\begin{Theorem}
Let $G$ be a three-dimensional unimodular Lie group with left-invariant Lorentzian metric. Then there exists a pseudo-orthonormal frame field $\left\{e_1, e_2, e_3\right\}$, such that the metric Lie algebra of $G$ is one of the following:
\begin{enumerate}
\item
\begin{equation*}
\mathcal{A}_1:
\begin{split}
[e_1,e_2] &= \lambda_3e_3, \\
[e_1,e_3] &= -\lambda_2e_2, \\
[e_2,e_3] &= \lambda_1e_1,
\end{split}
\end{equation*}
with $e_1$ timelike;
\item
\begin{equation*}
\mathcal{A}_2:
\begin{split}
[e_1,e_2] &= \left(1-\lambda_2\right)e_3-e_2, \\ 
[e_1,e_3] &= e_3-\left(1+\lambda_2\right)e_2, \\
[e_2,e_3] &= \lambda_1e_1,
\end{split}
\end{equation*}
with $e_3$ timelike;
\item
\begin{equation*}
\mathcal{A}_3:
\begin{split}
[e_1,e_2] &= e_1-\lambda e_3, \\ 
[e_1,e_3] &= -\lambda e_2-e_1, \\
[e_2,e_3] &= \lambda_1e_1+e_2+e_3,
\end{split}
\end{equation*}
with $e_3$ timelike;
\item
\begin{equation*}
\mathcal{A}_4:
\begin{split}
[e_1,e_2] &= \lambda_3e_2, \\ 
[e_1,e_3] &= -\beta e_1-\alpha e_2, \\
[e_2,e_3] &= -\alpha e_1+\beta e_2,
\end{split}
\end{equation*}
with $e_1$ timelike and $\beta \ne 0$.
\end{enumerate}
\end{Theorem}

\begin{Remark}
There are exactly six nonisomorphic three-dimensional unimodular Lie algebras and the corresponding types of three-dimensional unimodular Lie groups (see~\cite{Milnor}). All of them are listed in the Table~\ref{tab_PastukhovaSV:Chibrikova2006} together with conditions on structure constants for~which the Lie algebra has this type. If there is a ``$-$'' in the Table~\ref{tab_PastukhovaSV:Chibrikova2006} at the intersection of the row, corresponding to the Lie algebra, and the column, corresponding to the type, then it means that this type of the basis  is impossible for given Lie algebra. For~the~case of Lie algebra $\mathcal{A}_1$ we give only the signs of the triple $\left(\lambda_1,\lambda_2,\lambda_3\right)$ up to reorder and sign change.
\end{Remark}

\begin{table*}[ht!]
\caption{Three-dimensional unimodular Lie algebras}
\label{tab_PastukhovaSV:Chibrikova2006}
\begin{center} 
\begin{tabular}[c]{|c|c|c|c|c|} 
  \hline
  \multirow{2}*{Lie algebra} & \multicolumn{4}{c|}{Restrictions on the structure constants}  \\
   
  \cline{2-5}   & $\mathcal{A}_1$ & $\mathcal{A}_2$ & $\mathcal{A}_3$ & $\mathcal{A}_4$ \\
\hline
  $su(2)$            & $\left(+,+,+\right)$ & $-$ & $-$ & $-$ \\
\hline
  $sl(2,\mathbb{R})$ & $\left(+,+,-\right)$ & $\lambda_1 \ne 0$, $\lambda_2 \ne 0$ & $\lambda \ne 0$ & $\lambda_3 \ne 0$ \\
\hline
  $e(2)$             & $\left(+,+,0\right)$ & $-$ & $-$ & $-$ \\
\hline
  $e(1,1)$           & $\left(+,-,0\right)$ & \parbox[c]{80pt}{\centering $\lambda_1 = 0$, $\lambda_2 \ne 0$ \\ or \\ $\lambda_1 \ne 0$, $\lambda_2 = 0$} & $\lambda = 0$ & $\lambda_3 = 0$ \\
\hline
  $h$                & $\left(+,0,0\right)$ & $\lambda_1 = 0$, $\lambda_2 = 0$ & $-$ & $-$  \\
\hline
  $\mathbb{R}^{3}$   & $\left(0,0,0\right)$ & $-$ & $-$ & $-$  \\

  \hline
\end{tabular}
\end{center}
\end{table*}

\begin{Remark}
We note that similar bases was also constructed by G.~Calvaruso, L.\,A.~Cordero and P.\,E.~Parker~in~\cite{KlepikovaSV:Calvaruso2007,CP1997}.
\end{Remark}


\section{The isotropy of Schouten\nobreakdash--Weyl tensor}

In the paper \cite{R-S-Sh} the problem of existence of pseudo\nobreakdash-Rieman\-nian metrics, for which the squared length of the Schouten\nobreakdash--Weyl tensor is zero and some of the components of the~Schouten\nobreakdash--Weyl tensor are not zero sumultaneously has been tasked. This problem has been solved for the cases $\mathcal{A}_1$ and $\mathcal{A}_4$, and for nonunimodular three-dimensional Lie~algebras.
This paper presents a solutions to the problem for the cases $\mathcal{A}_2$ and $\mathcal{A}_3$.
 
\begin{Theorem} \label{Theorem:2}
Let $G$ be an unimodular three-dimensional Lie group with left-invariant Lorentzian metric, Lie algebra $\mathcal{A}_2$ or $\mathcal{A}_3$ and isotropic Schouten\nobreakdash--Weyl tensor. Then the~Lie~algebra of group~$G$ is~isomorphic to $e(1,1)$ or $sl(2,\mathbb{R})$, and the restrictions on~structural constants are contained in~the~Table~\ref{table:2}.
\end{Theorem}

\begin{table}[!ht]
\caption{Three-dimensional unimodular metric Lie algebras~$\mathcal{A}_2$ and $\mathcal{A}_3$ with the isotropic Schouten\nobreakdash--Weyl tensor}
\label{table:2}
\begin{center}
\begin{tabular}[c]{|c|c|c|}
\hline
Lie algebra & Restrictions on the structure constants & Type of Lie algebra\\
\hline 

\multirow{2}{*}{$\mathcal{A}_2$} & $\lambda_1=0$, $\lambda_{2}\ne 0$ & $e(1,1)$ \\
\cline{2-3}
                                 & $\lambda_1=\lambda_2\ne0$         & \multirow{2}{*}{$sl(2,\mathbb{R})$} \\
\cline{1-2}

$\mathcal{A}_3$                  & $\lambda\neq 0$                   & \\
\hline
\end{tabular}
\end{center}
\end{table}

\begin{proof}
We consider the case $\mathcal{A}_2$.  
Let be $\{E_1, E_2, E_3\}$ a basis, which was listed in the Theorem~\ref{Theorem:2}. 
Calculating the components of the Schouten\nobreakdash--Weyl tensor with help of the previously presented formulas, we see, that the non-trivial components of the Schouten\nobreakdash--Weyl tensor have the form
 \begin{equation}\label{sw1}
\begin{split}
&SW_{132}=-\lambda_1^3+\lambda_1^2\lambda_2,\\
&SW_{221}=SW_{331}=-\frac{1}{2}\lambda_1^2-2\lambda_2\lambda_1+4\lambda_2^2,\\
&SW_{231}=-\frac{1}{2}\lambda_1^2-2\lambda_2\lambda_1+4\lambda_2^2-\frac{1}{2}\lambda_1^3+\frac{1}{2}\lambda_1^2\lambda_2,\\
&SW_{321}=-\frac{1}{2}\lambda_1^2-2\lambda_2\lambda_1+4\lambda_2^2+\frac{1}{2}\lambda_1^3-\frac{1}{2}\lambda_1^2\lambda_2,
\end{split}
\end{equation}
and the squared length of the Schouten\nobreakdash--Weyl tensor is equal to
$$\|SW\|^2=-3\lambda_1^4(\lambda_1-\lambda_2)^2.$$
The formula shows that equality to the zero is achieved if: $\lambda_1=0$ or $\lambda_1=\lambda_2$. 
Thus the~Schouten\nobreakdash--Weyl tensor is not trivial if and only if $\lambda_1=0, \lambda_2 \neq 0$ or $\lambda_1=\lambda_2 \neq 0$. 

Now we consider the case $\mathcal{A}_3$.
Let be $\{E_1, E_2, E_3\}$ is a basis, which was listed in the Theorem~\ref{Theorem:2}. 
As in the case of $\mathcal{A}_2$, we calculate the components of the Schouten\nobreakdash--Weyl tensor
\begin{equation}\label{sw2}
\begin{split}
&SW_{121}=-SW_{131}=SW_{232}=-SW_{332}=\frac{3}{2}\lambda^2,\\
&SW_{221}=-SW_{231}=-SW_{321}=SW_{331}=-6\lambda.
\end{split}
\end{equation}
Thus the squared length of the Schouten\nobreakdash--Weyl tensor is trivial for any value of the structure constant $\lambda$, and the Schouten\nobreakdash--Weyl tensor can be zero only if $\lambda\equiv0$. 
\end{proof}

\begin{Remark}
The case of Lie algebra $\mathcal{A}_1$ or $\mathcal{A}_4$ and case of non-unimodular Lie algebras has been considered~in~\cite{R-S-Sh}.
\end{Remark}

\section{Almost harmonicity of the Schouten\nobreakdash--Weyl tensor}

In this part we investigate the three-dimensional metric Lie groups with zero curl and divergence of the Schouten\nobreakdash--Weyl tensor.

\begin{Theorem} \label{theorem:3}
Let $G$ be an unimodular three-dimensional Lie group with left\nobreakdash-invariant Lorentzian metric and Lie algebra $\mathcal{A}_2$ or $\mathcal{A}_3$. Then $\pdiv (SW)\equiv 0$. If moreover ${\pcurl(SW)=0}$, i.e.~the~Schouten\nobreakdash--Weyl tensor is almost harmonic, then the Schouten\nobreakdash--Weyl tensor is~trivial (so~$\qdiv(SW)=0$), and the Lie algebra of group $G$ is contained in the Table~\ref{table:3}.
\end{Theorem}

\begin{table}[!ht]
\caption{The three-dimensional unimodular metric Lie algebras $\mathcal{A}_2$ or $\mathcal{A}_3$ with trivial Schouten\nobreakdash--Weyl tensor}
\label{table:3}
\begin{center}
\begin{tabular}[c]{|c|c|c|}
\hline
Lie algebra & Restrictions on the structure constants & Type of Lie algebra\\
\hline 

$\mathcal{A}_2$ & $\lambda_1=\lambda_2=0$ & $h$ \\
\hline

$\mathcal{A}_3$ & $\lambda= 0$            & $e(1,1)$ \\
\hline
\end{tabular}
\end{center}
\end{table}

\begin{proof}
We consider the case of Lie algebra $\mathcal{A}_2$. Direct calculations show that~${\pdiv (SW)\equiv 0}$. 

Further we assume that the Schouten\nobreakdash--Weyl tensor is almost harmonic, i.e.~${\pcurl(SW)=0}$, so we have the system of equations:
\begin{equation}\label{curl1}
\begin{gathered}
-2\lambda_2\lambda_1^2-\lambda_1^3+16\lambda_2^2\lambda_1-16\lambda_2^3=0, \\
\left(\lambda_1^2+4\lambda_1\lambda_2-8\lambda_2^2\right)\left(\lambda_1-2\lambda_2\right)=0, \\
-\lambda_1^3+2\lambda_2\lambda_1^2+8\lambda_2^2\lambda_1+3\lambda_1^4-3\lambda_2\lambda_1^3=0, \\
-7\lambda_1^3+2\lambda_2\lambda_1^2+8\lambda_2^2\lambda_1=0, \\
7\lambda_1^3-2\lambda_2\lambda_1^2-3\lambda_1^4+3\lambda_2\lambda_1^3-8\lambda_2^2\lambda_1=0.
\end{gathered}
\end{equation}
Solving the system of equations (\ref{curl1}) for the structure constants of the Lie algebra $\mathcal{A}_2$, we obtain the solution ${\lambda_1=\lambda_2=0}$. 
It is immediately verified that for this solution the~Schouten\nobreakdash--Weyl tensor (\ref{sw1}) is trivial.



Next we consider the case of Lie algebra $\mathcal{A}_3$. Direct calculations show that in this case we have $\pdiv (SW)\equiv 0$. 

Let the Schouten\nobreakdash--Weyl tensor satisfies the condition $\pcurl(SW)=0$. In this case, we obtain the following system of equations:
\begin{equation}\label{curl2} 
\begin{gathered}
\lambda=0, \quad \lambda^2=0, \quad \lambda^3=0, \\
16\lambda+\lambda^3=0, \quad 16\lambda^2-\lambda^3=0, \quad \lambda^3-6\lambda=0, \quad \lambda^3+6\lambda=0.
\end{gathered}
\end{equation}
Solving the system of equations (\ref{curl2})  for the structure constants of the Lie algebra $\mathcal{A}_3$, we obtain the solution $\lambda=0$.
It is easy to verify, that for this solution the~Schouten\nobreakdash--Weyl tensor (\ref{sw2}) is trivial.

\end{proof}

\begin{Remark}
The case of Lie algebra $\mathcal{A}_1$ or $\mathcal{A}_4$ and case of non-unimodular Lie algebras has been considered~in~\cite{G-R-S-2}.
\end{Remark}

\section{The harmonicity of contraction of the Schouten\nobreakdash--Weyl tensor}

In this section we consider the question about the harmonicity of contraction of the Schouten\nobreakdash--Weyl tensor, which was defined by~(\ref{ev0}), through the Schouten\nobreakdash--Weyl tensor $SW$ and some vector field $V$.

\begin{Theorem} \label{theorem:4}
Let $G$ be an unimodular three-dimensional Lie group with left\nobreakdash-invariant Lorentzian metric, Lie algebra $\mathcal{A}_2$ or $\mathcal{A}_3$ and $w$ is a harmonic tensor. Then the restrictions on structure constants of the Lie algebra of group $G$ and components of vector field $\{V^{k}\}$ are contained in the Table~\ref{table:4}. 
\end{Theorem}

\begin{table}[!ht]
\caption{The three-dimensional unimodular metric Lie algebras $\mathcal{A}_2$ or $\mathcal{A}_3$ with harmonic contraction of the Schouten\nobreakdash--Weyl tensor}
\label{table:4}
\begin{center} \small
\begin{tabular}[c]{|c|c|c|}
\hline
\parbox{90pt}{\centering Restrictions on the structure constants} & \parbox{65pt}{\centering Type of\\ Lie algebra} &  Vector field\\
\hline 

\multicolumn{3}{|c|}{Lie algebra $\mathcal{A}_2$}\\
\hline
$\lambda_1=-2$, $\lambda_2=0$ & $e(1,1)$ & $V=(V^1,V^2,V^3)$, $V^i \in \mathbb{R}$  is harmonic. \\
\hline
$\lambda_{i}=L_i(V^2,V^3)$ & \parbox{65pt}{\centering $sl(2,\mathbb{R})$\\ or $e(1,1)$}  & \parbox{190pt}{\centering $V=(0,V^2,V^3)$, $V^2,V^3 \in \mathbb{R}\setminus\{0\}$,\\ $V^2 \neq \left(-7\pm4\sqrt{3}+2\sqrt{24+14\sqrt{3}}\right)V^3$\\ is harmonic.} \\
\hline
$\lambda_2=0$, $\lambda_{1}\neq 0$ & $e(1,1)$ & $V=(0,V^2,V^2)$, $V^2 \in \mathbb{R}$  is harmonic. \\
\hline
$\lambda_1=0$, $\lambda_{2}\neq 0$ & $e(1,1)$ & $V=(V^1,0,0)$, $V^1 \in \mathbb{R}$ \\
\hline 

\multicolumn{3}{|c|}{Lie algebra $\mathcal{A}_3$}\\
\hline
$\lambda=0$ & $e(1,1)$ & $V=(0,-V^3,V^3)$, $V^3 \in \mathbb{R}$ is harmonic. \\
\hline
$\lambda=L_3$ & $sl(2,\mathbb{R})$ & \parbox{190pt}{\centering $V=(V^1,V^2,V^3)$,  $V^3 \in \mathbb{R}$,\\ $V^1=f(\lambda)V^3$, $V^2=h(\lambda)V^3$} \\
\hline
\end{tabular}
\end{center}
{\footnotesize where 
\begin{gather*}
L_1(x,y)=\frac{1}{6}\left(F\left(x,y\right)+\frac{H\left(x,y\right)}{F\left(x,y\right)}-\left(x+y\right)^2\right),  \quad L_2\left(x,y\right)=\frac{y^2-x^2}{2xy}, \\
F\left(x,y\right)=\left(P\left(x,y\right)+6\sqrt{6\left(x^2+y^2\right)\left(x-y\right)^4Q\left(x,y\right)}\right)^{1/3}, \\
P\left(x,y\right)=53x^6-150x^5y+507x^4y^2-308x^3y^3 + 507x^2y^4-150xy^5+53y^6, \\
Q\left(x,y\right)=13x^6-22x^5y+163x^4y^2-52x^3y^3 + 163x^2y^4-22xy^5+13y^6, \\
H\left(x,y\right)=x^4+28x^3y+6x^2y^2+28xy^3+y^4, \\
L_3=\pm\frac{1}{4}\sqrt{\sqrt{6\left(136643+512A+\frac{104960}{A}+\frac{85821417}{B}\right)}-483-B} \approx \pm 89.072, \\
A=\left(3763+6\sqrt{154029}\right)^{1/3}\approx 18.289, \\
B=\frac{1}{A}\sqrt{409929A-3072A^2-629760}\approx 565.076, \\
f\left(\lambda\right)=\frac{1}{8192}\lambda\left(440\lambda^6+9047\lambda^4-47248\lambda^2+30400\right), \\
h\left(\lambda\right)=-\frac{1}{2048}\left(88\lambda^6+2219\lambda^4+4322\lambda^2-2112\right). 
\end{gather*} }
\end{table}

\begin{proof}
First we consider the case of Lie algebra $\mathcal{A}_2$. 
Using (\ref{ev0}), we find the nontrivial components of the~tensor $w_{ij}$
\begin{gather*}
\begin{split}w_{12}=\frac{1}{2}V^2\lambda_1^2-2V^2\lambda_1&-6V^2-2V^2\lambda_2^2+2V^3+ 2V^3\lambda_2^2-2V^3\lambda_2\lambda_1 +\\
&+\frac{1}{2}V^3\lambda_2\lambda_1^2+ \frac{1}{2}V^3\lambda_1^3+\frac{1}{2}V^3\lambda_1^2- 6V^3\lambda_2-  2V^3\lambda_2^3 ,\end{split}\\
\begin{split}w_{13}= -2V^2\lambda_2^2-2V^2\lambda_2\lambda_1&+\frac{1}{2}V^2\lambda_2\lambda_1^2-\frac{1}{2}V^2\lambda_1^3-\frac{1}{2}V^2\lambda_1^2-6V^2\lambda_2-\\
&-2V^2\lambda_2^3-2V^2- \frac{1}{2}V^3\lambda_1^2+2V^3\lambda_1+6V^3+2V^3\lambda_2^2,\end{split}\\
w_{23}=-V^1(\lambda_1^3+\lambda_1^2+4\lambda_2^2+4).
\end{gather*}
Computing the curl and divergence of the tensor $w_{ij}$, we see that the equalities
$\mathrm{curl}(w)\equiv 0$ and \hbox{$\mathrm{div}(w)\equiv0$} are equivalent to the system of equations
\begin{equation}\label{lev2s}
\begin{gathered}
-V^1(\lambda_1^3+\lambda_1^2+4\lambda_2^2+4)\lambda_1 =0, \\
\begin{split}4V^2\lambda_2\lambda_1+V^2\lambda_2\lambda_1^3-V^2\lambda_2^2\lambda_1^2+4V^2\lambda_2^2\lambda_1+16V^2-16V^3+20V^2\lambda_2^2-8V^3\lambda_2^2-\\
-V^3\lambda_1^3+V^2\lambda_1^3+8V^2\lambda_2^3+4V^2\lambda_2^4-4V^3\lambda_1+16V^2\lambda_2+4V^2\lambda_1=0,\end{split} \\
\begin{split}-4V^3\lambda_2\lambda_1+4V^3\lambda_2^2\lambda_1-V^3\lambda_2^2\lambda_1^2-V^3\lambda_2\lambda_1^3-16V^2+16V^3-8V^2\lambda_2^2+20V^3\lambda_2^2+\\
+V^3\lambda_1^3-8V^3\lambda_2^3-V^2\lambda_1^3+4V^3\lambda_2^4-16V^3\lambda_2+4V^3\lambda_1-4V^2\lambda_1=0.\end{split}
\end{gathered}
\end{equation}

Solving the system of equations (\ref{lev2s}) for the structure constants of the Lie algebra and the components of vector field $\{V^{k}\}$, we define Lie algebras and the corresponding directions, for which the tensor $w_{ij}$ is harmonic. We obtain the following solutions:

\begin{enumerate}
\item $V=(V^1,V^2,V^3)$, $V^i \in \mathbb{R}$, $\lambda_1=-2$, $\lambda_2=0$;
\item ${V=(0,V^2,V^3)}$, ${ V^2,V^3 \in \mathbb{R}\setminus\{0\}}$, 
${\lambda_1=\frac{1}{6}\left( F(V^2,V^3)+\frac{H(V^2,V^3)}{F(V^2,V^3)}-(V^2+V^3)^2 \right)}$,  
${\lambda_2=\frac{(V^3)^2-(V^2)^2}{2V^2V^3}}$,
${V^2 \neq \left(-7\pm4\sqrt{3}+2\sqrt{24+14\sqrt{3}}\right)V^3}$, where
{\small \begin{gather*}
F\left(x,y\right)=\left(P\left(x,y\right)+6\sqrt{6\left(x^2+y^2\right)\left(x-y\right)^4Q\left(x,y\right)}\right)^{1/3}, \\
P\left(x,y\right)=53x^6-150x^5y+507x^4y^2-308x^3y^3 + 507x^2y^4-150xy^5+53y^6, \\
Q\left(x,y\right)=13x^6-22x^5y+163x^4y^2-52x^3y^3 + 163x^2y^4-22xy^5+13y^6, \\
H\left(x,y\right)=x^4+28x^3y+6x^2y^2+28xy^3+y^4;
\end{gather*} }
\item $V=(0,V^2,V^2)$, $V^2 \in \mathbb{R}$, $\lambda_1 \in \mathbb{R}\setminus\{0\}$, $\lambda_2=0$;
\item $ V=(V^1,0,0)$, $\lambda_1=0$,  $\lambda_2 \in \mathbb{R}\setminus\{0\}$, $V^1 \in \mathbb{R}$.
\end{enumerate}

Let $\{V^{k}\}$ be harmonic. Computing the curl and divergence of the vector field $\{V^{k}\}$, we see that the equalities $\mathrm{curl}(V)\equiv 0$ and \hbox{$\mathrm{div}(V)\equiv0$} are equivalent to the equation:
\begin{equation}\label{divVcurlV_2}
V^1(2+\lambda_1)=0.
\end{equation}
Solving the system of equations (\ref{lev2s}) and (\ref{divVcurlV_2}), we obtain, that the directions~1\nobreakdash--3 are harmonic and the direction 4 isn't harmonic. 

Now we consider the case of Lie algebra $\mathcal{A}_3$. 
Using (\ref{ev0}), we find the nontrivial components of the tensor $w_{ij}$
\begin{gather*}
w_{12}=-\frac{11}{2}V^1\lambda^2+4V^1+2V^2\lambda-\lambda^3V^3+2V^3\lambda ,\\
w_{13}=\frac{3}{2}V^1\lambda^2-4V^1-2\lambda^3V^2+6V^2\lambda-2V^3\lambda,\\
w_{23}=-\lambda^3V^1+4V^1\lambda+\frac{3}{2}V^2\lambda^2-4V^2-4V^3+\frac{11}{2}V^3\lambda^2 .
\end{gather*}
Computing the curl and divergence of the tensor $w_{ij}$, we see that the equalities $\mathrm{curl}(w)\equiv 0$ and \hbox{$\mathrm{div}(w)\equiv0$} are equivalent to the system of equations:
\begin{equation}\label{lev3}
\begin{gathered}
-16V^1+22V^1\lambda^2-\lambda^3V^2-16V^3\lambda-2\lambda^4V^1+13\lambda^3V^3=0, \\
-\lambda^3V^1+4\lambda^4V^2-15V^2\lambda^2-7V^3\lambda^2+8V^2+8V^3=0, \\
-16V^1\lambda+13\lambda^3V^1-7V^2\lambda^2+2\lambda^4V^3-15V^3\lambda^2+8V^2+8V^3=0.
\end{gathered}
\end{equation}

Solving this system of equations for the structure constants of the Lie algebra and the~components of vector field $\{V^{k}\}$, we obtain the following solutions:
\begin{enumerate}
\item $V=(0, -V^3, V^3)$, $V^3 \in \mathbb{R}$,  $\lambda=0$;
\item $V=(V^1,V^2,V^3)$, $V^3\ne0$ and 
{\footnotesize\begin{gather*}
V^1=\frac{1}{8192}V^3\lambda(440\lambda^6+9047\lambda^4-47248\lambda^2+30400),\\
V^2=-\frac{1}{2048}V^3(88\lambda^6+2219\lambda^4+4322\lambda^2-2112),\\
\lambda=\pm\frac{1}{4}\sqrt{\sqrt{6\left(136643+512A+\frac{104960}{A}+\frac{85821417}{B}\right)}-483-B} \approx \pm 89.072;\\
A=(3763+6\sqrt{154029})^{1/3}\approx 18.289,\\
B=\frac{1}{A}\sqrt{409929A-3072A^2-629760}\approx 565.076.
\end{gather*} }
\end{enumerate}

Let $\{V^{k}\}$ be harmonic. Computing the curl and divergence of the vector field $\{V^{k}\}$, we see that the equalities $\mathrm{curl}(V)\equiv 0$ and \hbox{$\mathrm{div}(V)\equiv0$} are equivalent to the equation:
\begin{equation}\label{divVcurlV_3}
-\lambda V^1+2V^2+2V^3=0.
\end{equation}
Solving the system of equations (\ref{lev3}) and (\ref{divVcurlV_3}), we obtain, that the direction~1 is harmonic and the direction 2 isn't harmonic. 
\end{proof}

\begin{Remark}
A similar theorems for the contraction of the Schouren-Weyl tensor in~the~cases $\mathcal{A}_1$ and $\mathcal{A}_4$ and also in the case of non-unimodular Lie~algebras has been proved in~\cite{G-R-S-2}.
\end{Remark}


\Addresses

\end{document}